\documentclass[11pt]{amsart}
\usepackage{amsmath,amssymb,mathrsfs}
\newcommand{\doi}[1]{\href{https://doi.org/#1}{\nolinkurl{doi:#1}}}
\usepackage[hidelinks]{hyperref}
\newtheorem{theorem}{{Theorem}}
\newtheorem{lemma}[theorem]{Lemma}
\newtheorem{proposition}[theorem]{Proposition}
\theoremstyle{remark}
\newtheorem*{remark}{Remark}
\newcommand{\tr}{\on{{tr}}}
\newcommand{\la}{\langle}
\newcommand{\rg}{\rangle}

\newcommand{\bR}{\mathbb{R}} 
\newcommand{\nc}{\newcommand}
\nc{\on}{\operatorname}
\nc{\vol}{\text{\rm vol}}

\newcommand{\gd}{\nabla}

\newcommand{\chib}{\underline{\chi}}
\newcommand{\trch}{\tr(\chi)}
\newcommand{\trchb}{\tr(\chib)}
\newcommand{\detchb}{\det(\chib)}
\newcommand{\Lb}{\underline L}

\title{The Gibbons--Penrose inequality in Minkowski spacetime}
\author{Jingbo Wan}
 \address{Laboratoire Jacques-Louis Lions de Sorbonne Universit\'e, 4 place Jussieu, Paris 75005, France} \email{jingbo.wan@sorbonne-universite.fr}

\author{Yipeng Wang}
\address{Department of Mathematics, Princeton University, Princeton, NJ 08544, USA}
\email{yipeng.wang@princeton.edu}
\thanks{J.W. is supported by ERC-2023 AdG 101141855 BlaHSt.}
\usepackage{xcolor}

\makeatletter
\let\GPblacktagform\tagform@
\renewcommand{\tagform@}[1]{\textcolor{black}{\GPblacktagform{#1}}}
\makeatother
\begin{document}

\begin{abstract}
    We prove the Gibbons--Penrose inequality for smooth closed spacelike codimension-two submanifolds of Minkowski spacetime whose past null second fundamental form is non-negative. The proof is motivated by the Alexandrov--Bakelman--Pucci technique.
\end{abstract}
\maketitle
\section{Introduction}
Penrose~\cite{Penrose} proposed an energy--area inequality for collapsing null dust shells, motivated by cosmic censorship. For a shell with flat interior, the matter density can be prescribed so that a spacelike section is marginally outer trapped with respect to the exterior geometry. The proposed bound then reduces to a geometric inequality for the section in Minkowski spacetime, known as the Gibbons--Penrose inequality; see~\cite{BW,Gibbons,MS12}. We prove a strengthened version in every spacetime dimension at least four, assuming that the past null second fundamental form is nonnegative.

Let $n\geq 3$ be an integer, and let $\Sigma^{n-1}\subset \bR^{n,1}$ be a closed spacelike submanifold with induced Riemannian metric $g$. We consider a fixed future-directed unit timelike vector $T_0$ with $\la T_0,T_0\rg=-1$. Let $\vec H$ denote the mean curvature vector of $\Sigma$, with the convention that the mean curvature vector of a standard round sphere is inward pointing. Let $L$ and $\Lb$ be the two null normal vector fields along $\Sigma$ satisfying $\la L,\Lb\rg=2$, where $L$ is future-directed and $\Lb$ is past-directed. We normalize the null frame by requiring $\la \Lb,T_0\rg=1$. In particular,
\[
\vec H = \frac12 ( \la\vec H,L\rg\,\Lb + \la\vec H,\Lb\rg\,L ).
\]
The dual mean curvature vector $\vec J$ is defined by
\[
\vec J = \frac12 ( \la\vec H,L\rg\,\Lb - \la\vec H,\Lb\rg\,L ).
\]
We define the null second fundamental forms as follows. Given a local orthonormal frame $\{e_1,\ldots,e_{n-1}\}$ of $\Sigma$, we set $ \chi(e_i,e_j) = \la\bar D_{e_i}L,e_j\rg $ and $\chib(e_i,e_j) = \la\bar D_{e_i}\Lb,e_j\rg$, where $\bar D$ denotes the Levi-Civita connection of $\bR^{n,1}$. Consequently, we have $\trch=-\la\vec H,L\rg$ and $\trchb=-\la\vec H,\Lb\rg$. 
\begin{theorem}\label{thm:GP}
Suppose that $\chib\geq 0$ at each point of $\Sigma$. Then
\begin{align}\label{eq:strengthened-GP}
-\int_\Sigma \la\vec J,T_0\rg\,d\mu &\geq (n-1)|S^{n-1}|^{\frac1{n-1}} |\Sigma|^{\frac{n-2}{n-1}} \notag\\
&\quad +\frac14\int_\Sigma \trchb\,|T_0^{\tan}|^2\,d\mu +\frac34\int_\Sigma \chib(T_0^{\tan},T_0^{\tan})\,d\mu.
\end{align}
In particular,
\begin{equation}\label{eq:GP}
-\int_\Sigma \la\vec J,T_0\rg\,d\mu \geq (n-1)|S^{n-1}|^{\frac1{n-1}} |\Sigma|^{\frac{n-2}{n-1}}.
\end{equation}
Moreover, equality holds in \eqref{eq:GP} if and only if $\Sigma$ is a round sphere contained in a spacelike hyperplane orthogonal to $T_0$.
\end{theorem}
\begin{remark}
The original Gibbons--Penrose conjecture assumes that $\Sigma$ is \emph{past null convex}: the null generators starting from $\Sigma$ in the direction $\Lb$ form a smooth null hypersurface extending to past null infinity; see~\cite[Conjecture~1]{BW}. This condition implies $\chib\geq0$, so Theorem~\ref{thm:GP} applies to this class.
\end{remark}
The normalization $\la\Lb,T_0\rg=1$ and the identity $\int_\Sigma\vec H\,d\mu=0$ give
\[
-\int_\Sigma\la\vec J,T_0\rg\,d\mu =\int_\Sigma\trch\,d\mu.
\]
For $n=3$, the preceding identity rewrites ~\eqref{eq:GP} as
\[
\int_\Sigma\trch\,d\mu\geq\sqrt{16\pi|\Sigma|}.
\]
If $\Sigma$ lies in a spacelike hyperplane orthogonal to $T_0$, the theorem reduces to the classical Minkowski inequality for convex hypersurfaces.

Gibbons~\cite{Gibbons} related the collapsing-shell inequality to the classical Minkowski inequality and proposed a proof by orthogonal projection onto a spacelike hyperplane. The projected mean curvature calculation, however, misses a term involving the normal connection; see~\cite{BW,MS12}. Tod~\cite{Tod} proved the inequality for surfaces in the past null cone of a point by a Sobolev inequality. Mars and Soria~\cite{MS12} analyzed the gap in the projection argument, expressed the inequality using a convex constant-time section of the incoming null hypersurface and a time-height function, and extended the null-cone case to arbitrary dimension using Beckner's sharp inequality on the sphere{~\cite{Beckner}}. Brendle and Wang~\cite{BW} formulated the inequality in terms of the dual mean curvature vector $\vec J$. Moreover, they proposed an extension of the Gibbons--Penrose conjecture to Schwarzschild spacetime, and proved several cases, including surfaces with strictly convex spatial projection in the Minkowski limit. Mars and Soria~\cite{MS13} obtained a related sufficient condition using normal graphs in Euclidean space.

Our proof uses the Alexandrov--Bakelman--Pucci technique. Cabr\'e~\cite{Cabre} applied this method to the classical isoperimetric inequality by solving a linear Neumann problem and estimating the Jacobian determinant of its gradient map on a contact set. Cabr\'e, Ros-Oton, and Serra~\cite{CRS} used this method to prove sharp weighted and anisotropic isoperimetric inequalities. Brendle~\cite{Brendle} adapted the method to submanifolds and obtained a sharp Michael--Simon Sobolev inequality in codimension at most two. Tsai and Wang~\cite{TW} proved related inequalities for spacelike submanifolds in Minkowski spacetime. For a recent survey of the method and its applications to geometric inequalities, see Brendle~\cite{BrendleSurvey}.

After approximation, we may assume that $\chib>0$. The normalized null Gauss map is then a diffeomorphism onto the round sphere on each connected component. We then solve a linear equation in divergence form on $\Sigma$ and transport its solution to the sphere. Using this solution, we construct a map $\Phi$ that sends its contact set onto $S^{n-1}$ and control its Jacobian determinant there.
\subsection*{Disclosure of AI tools.}
The authors conceived the application of the Alexandrov--Bakelman--Pucci technique and developed the mathematical arguments presented in this paper. ChatGPT 6 Pro was used to assist with calculations, the analysis of proof details, and revision of the exposition. The authors take full responsibility for the mathematical content and the final manuscript.
\section{Proof of Theorem \ref{thm:GP}}
We first describe an auxiliary inequality.
\begin{proposition}\label{prop:aux-inequality}
Suppose that $\chib>0$, then
\begin{align*}
&\int_{\Sigma} \Big( \trchb +\frac34\trchb\,|T_0^{\tan}|^2- \frac34 {\chib} (T_0^{\tan},T_0^{\tan}) \Big)\,d\mu\notag\\
&\geq (n-1)|S^{n-1}|^{\frac1{n-1}} |\Sigma|^{\frac{n-2}{n-1}}.
\end{align*}
\end{proposition}
We first show that Proposition \ref{prop:aux-inequality} implies Theorem \ref{thm:GP}.
\begin{proof}[Proof of Theorem \ref{thm:GP}]
Given $\Sigma$ with $\chib\geq0$, we may choose a sequence of spacelike submanifolds $\Sigma_j\to\Sigma$ smoothly with $\chib_{\Sigma_j}>0$. Hence, it suffices to prove the inequalities when $\chib>0$, which we assume for now.

Since $\la\Lb,T_0\rg=1$, we may write
\begin{align*}
T_0 = T_0^{\tan} +\frac12L -\frac12 \big( 1+|T_0^{\tan}|^2 \big)\Lb.
\end{align*}
On the other hand, since $T_0$ is a parallel vector field in $\bR^{n,1}$, we have
\[
0 = \operatorname{div}_\Sigma (T_0^{\tan}) +\frac12\trch -\frac12 \big( 1+|T_0^{\tan}|^2 \big) \trchb.
\]
Combining these identities gives
\begin{align*}
-\la\vec J,T_0\rg &= \frac12\trch + \frac12 \big( 1+|T_0^{\tan}|^2 \big) \trchb\\
&=\big( 1+|T_0^{\tan}|^2 \big) \trchb-\operatorname{div}_\Sigma (T_0^{\tan}).
\end{align*}
Integrating over $\Sigma$ and using the divergence theorem gives
\begin{align*}
&-\int_\Sigma \la\vec J,T_0\rg\,d\mu\\ &= \int_\Sigma \big( 1+|T_0^{\tan}|^2 \big) \trchb\,d\mu\\
&=\int_{\Sigma} \Big( \trchb +\frac34\trchb\,|T_0^{\tan}|^2- \frac34 {\chib} (T_0^{\tan},T_0^{\tan}) \Big)\,d\mu\\
&\qquad + \frac14\int_{\Sigma}\trchb\,|T_0^{\tan}|^2\,d\mu + \frac34\int_{\Sigma} \chib(T_0^{\tan},T_0^{\tan})\,d\mu\\
&\geq (n-1)|S^{n-1}|^{\frac1{n-1}} |\Sigma|^{\frac{n-2}{n-1}}\\
&\qquad + \frac14\int_{\Sigma}\trchb\,|T_0^{\tan}|^2\,d\mu + \frac34\int_{\Sigma} \chib(T_0^{\tan},T_0^{\tan})\,d\mu.
\end{align*}
Here, the last inequality follows from Proposition \ref{prop:aux-inequality}. Passing to the limit for proves both inequalities for $\chib\geq0$.

Suppose now that equality holds in~\eqref{eq:GP}. Applying the inequality to each connected component shows that $\Sigma$ is connected. Moreover, we have $\trchb\,|T_0^{\tan}|^2=0$ at each point of $\Sigma$. We may rescale $\Sigma$ such that $|\Sigma|=|S^{n-1}|$, and we have $\int_{\Sigma}\trchb\,d\mu=(n-1)|\Sigma|$. 

Next, we fix any $\varphi\in C^{\infty}(\Sigma)$ and a real number $t$. We define $\Sigma_t=\{x+t\varphi(x)\underline{L}:x\in\Sigma\}$. If $|t|$ is sufficiently small, we have $\chib_{\Sigma_t}\geq 0$ and hence \eqref{eq:GP} applies to each $\Sigma_t$. Differentiating with respect to $t$ gives
\begin{align*}
    0&=\frac{d}{dt}\big|_{t=0}\Big(\int_{\Sigma_t} \la\vec J_{\Sigma_t},T_0\rg\,d\mu_{\Sigma_t}+(n-1)|S^{n-1}|^{\frac{1}{n-1}}|\Sigma_t|^{\frac{n-2}{n-1}}\Big)\\
    &=\int_{\Sigma}\varphi\,(|\chib|^2-\trchb^2+(n-2)\,\trchb)\,d\mu.
\end{align*}
Since $\varphi$ is arbitrary, we must have 
\[
    |\chib|^2-\trchb^2+(n-2)\,\trchb=0
\]
at each point of $\Sigma$. Since $|\chib|^2\geq \frac{\trchb^2}{n-1}$, we conclude that either $\trchb=0$ or $\trchb\geq n-1$. Since $\Sigma$ is connected and $\int_{\Sigma}\trchb\,d\mu=(n-1)|\Sigma|$, we must have $\trchb=n-1$ at each point of $\Sigma$. Using the preceding identity again, we know that $|\chib|^2=n-1$ at each point of $\Sigma$. Therefore, we conclude that $\chib=g$. Consequently, $T_0^{\tan}=0$ at each point of $\Sigma$. From this, we can easily deduce that $\Sigma$ is a round unit sphere in a spacelike hyperplane orthogonal to $T_0$. This completes the proof of the rigidity statement.
\end{proof}
In the rest of this section, we will prove Proposition \ref{prop:aux-inequality}. Without loss of generality, we assume that $\Sigma$ is connected. By scaling, we may further assume that
\begin{align}\label{eq:normalization}
\int_{\Sigma} \Big( \trchb +\frac34\trchb|T_0^{\tan}|^2- \frac34 {\chib} (T_0^{\tan},T_0^{\tan}) \Big)d\mu = (n-1)|\Sigma|.
\end{align}
Since $\chib>0$ and $\Sigma$ is connected, we can find a smooth function $u:\Sigma\to\bR$ such that
\begin{align}\label{eq:aux-PDE}
\operatorname{div}_{\Sigma} ( \chib^{-1}(\gd^\Sigma u) ) &= (n-1) -\trchb-\frac{3}{4}\trchb\,|T_0^{\tan}|^2 +\frac34 {\chib} (T_0^{\tan},T_0^{\tan}).
\end{align}
Next, we identify the unit sphere $S^{n-1}$ with
\[
S^{n-1} = \{ \xi\in\bR^{n,1}: \la\xi,\xi\rg=0,\ \la\xi,T_0\rg=1 \}.
\]
The metric induced from $\bR^{n,1}$ on this set is the standard round metric of sectional curvature $1$. By our normalization, the null Gauss map $ \Lb:\Sigma\mapsto S^{n-1} $ is a local diffeomorphism with Jacobian determinant $\detchb$. Since $\Sigma$ is compact and connected and $S^{n-1}$ is simply connected, it is a diffeomorphism. We may therefore define a smooth function $ \bar u:S^{n-1}\to\bR $ by $ \bar u(\Lb(p))=u(p). $ We now define
\[
U = \{ x\in S^{n-1}: |\gd^{S^{n-1}}\bar u(x)|<\pi \}
\]
and set $\Phi:U\to S^{n-1}$ by $\Phi(x) = \exp_x ( \gd^{S^{n-1}}\bar u(x) )$. We define $A$ to be the set of all $x\in U$ such that
\[
\bar u(y) +\frac12d_{S^{n-1}}(y,\Phi(x))^2 \geq \bar u(x) +\frac12|\gd^{S^{n-1}}\bar u(x)|^2
\]
for every $y\in S^{n-1}$.

\begin{lemma}\label{lem:surjective}
We have $ \Phi(A)=S^{n-1}. $
\end{lemma}
\begin{proof}
Fix $y\in S^{n-1}$, then there exists $\bar x\in S^{n-1}$ such that
\[
\bar u(\bar x)+\frac12d_{ S^{n-1}}(\bar x,y)^2 = \min_{x\in S^{n-1}} \Big( \bar u(x)+\frac12d_{ S^{n-1}}(x,y)^2 \Big).
\]
We claim that $\bar x\neq -2T_0-y$. Indeed, if $\bar x=-2T_0-y$, choose a unit vector $v\in T_{\bar x} S^{n-1}$ such that $\langle\nabla\bar u(\bar x),v\rangle\le0$. Along the geodesic $\gamma(t)=\exp_{\bar x}(tv)$, we have $ d_{ S^{n-1}}(\gamma(t),y)=\pi-t $ for $t>0$ sufficiently small, and hence
\[
\lim_{t\downarrow 0} \frac{1}{t}\Big( \bar u(\gamma(t)) +\frac12d_{ S^{n-1}}(\gamma(t),y)^2-\bar u(\bar x) -\frac{\pi^2}{2}\Big)=\la \gd^{S^{n-1}}\bar u(\bar x),v\rg-\pi <0.
\]
This is a contradiction.

Therefore, we have $d_{ S^{n-1}}(\bar x,y)<\pi$ and hence the function $d_{S^{n-1}}(y,\cdot)^2$ is smooth around $\bar x$. The first variation formula gives $ \nabla^{S^{n-1}}\bar u(\bar x)=\exp_{\bar x}^{-1}(y). $ Consequently, we have $\Phi(\bar x)=y$ and
\[
|\nabla^{S^{n-1}}\bar u(\bar x)| =d_{ S^{n-1}}(\bar x,y)<\pi.
\]
Moreover, by the minimizing property of $\bar x$, we have
\[
\bar u(x) +\frac12d_{ S^{n-1}}(x,\Phi(\bar x))^2 \ge \bar u(\bar x) +\frac12|\nabla\bar u(\bar x)|^2
\]
for every $x\in S^{n-1}$. Hence, $\bar x\in A$ and $\Phi(\bar x)=y$.
\end{proof}

\begin{lemma}\label{lem:determinant}
Let $\bar x\in A$, and set $ s=|\nabla^{ S^{n-1}}\bar u(\bar x)|. $ Then if $s>0$, we have
\[
D^2_{S^{n-1}}\bar u +s\,\frac{\cos (s)}{\sin(s)}\,g_{ S^{n-1}} +\Big(\frac{1}{s^2}-\frac{\cos(s)}{s\,\sin(s)}\Big)\, d\bar u\otimes d\bar u \geq 0
\]
at $\bar x$. Moreover, we have
\begin{align*}
&\det D\Phi(\bar x)= \Big(\frac{\sin s}{s}\Big)^{n-2}\notag\\
&\qquad \cdot\det\Big( D^2_{S^{n-1}}\bar u +s\,\frac{\cos (s)}{\sin(s)}\,g_{ S^{n-1}} +\Big(\frac{1}{s^2}-\frac{\cos(s)}{s\,\sin(s)}\Big)\, d\bar u\otimes d\bar u \Big).
\end{align*}
If $s=0$, then at $\bar x$ we have
$D^2_{S^{n-1}}\bar u+g_{S^{n-1}}\geq0$ and $\det D\Phi(\bar x)=\det(D^2_{S^{n-1}}\bar u+g_{S^{n-1}})$.
\end{lemma}

\begin{proof}
We fix $ y=\Phi(\bar x) =\exp_{\bar x}(\nabla^{S^{n-1}}\bar u(\bar x)). $ Since $\bar x\in A$, the function $ x\mapsto \bar u(x)+\frac12d_{ S^{n-1}}(x,y)^2 $ attains a minimum at $\bar x$. If $s=0$, then $y=\bar x$. We have $\frac12D^2_{S^{n-1}}d_{S^{n-1}}(\cdot,\bar x)^2=g_{S^{n-1}}$ at $\bar x$, so we have $0\leq D^2_{S^{n-1}}\bar u+g_{S^{n-1}}$. On the other hand, since $\langle D\Phi(\bar x)(X),Y\rangle=g_{S^{n-1}}(X,Y)+D^2_{S^{n-1}}\bar u(\bar x)(X,Y)$ for every $X,Y\in T_{\bar x}S^{n-1}$, both assertions follow. We may therefore assume that $s>0$.  Hence, $x\to d_{S^{n-1}}(x,y)^2$ is smooth at $x=\bar x$ and
\begin{align*}
0&\leq D^2_{S^{n-1}}\bar u(\bar x)+\frac{1}{2}D^2_{S^{n-1}}d_{S^{n-1}}(\bar x,y)^2\\
&=D^2_{S^{n-1}}\bar u +s\,\frac{\cos (s)}{\sin(s)}\,g_{ S^{n-1}} +\Big(\frac{1}{s^2}-\frac{\cos(s)}{s\,\sin(s)}\Big)\, d\bar u\otimes d\bar u.
\end{align*}
We next compute the Jacobian of $\Phi$ (cf.~\cite[Lemma~3.2]{CabreRiemannian}). We choose an orthonormal basis $\{e_1,\ldots,e_{n-1}\}$ of $T_{\bar x} S^{n-1}$ such that $ e_1=s^{-1}\nabla^{S^{n-1}}\bar u(\bar x). $ We then define $ \bar\gamma(t) = \exp_{\bar x}^{S^{n-1}}(t\,\nabla^{S^{n-1}}\bar u(\bar x)) $ for $0\leq t\leq 1$. For each $1\leq i\leq n-1$, let $E_i(t)$ be the unique parallel vector field along $\bar\gamma$ satisfying $ E_i(0)=e_i. $ Moreover, let $X_i(t)$ be the unique Jacobi field along $\bar\gamma$ satisfying $ X_i(0)=e_i$ and
\[
D_tX_i(0) = \sum_{j=1}^{n-1} (D^2_{S^{n-1}}\bar u)(e_i,e_j)\,e_j.
\]
Consequently, we have $ D\Phi(\bar x)(e_i)=X_i(1). $ Next, we  define the matrix $P(t)$ by $ P_{ij}(t)=\langle X_i(t),E_j(t)\rangle . $ It then follows that $\det D\Phi(\bar x)=\det (P_{ij}(1))$. Since the fields $E_j$ are parallel, we know that $P_{ij}(0)=\delta_{ij}$ and $P_{ij}'(0) = (D^2_{S^{n-1}}\bar u)(e_i,e_j)$. Moreover, the Jacobi equation gives $ P_{i1}''(t)=0 $ and
\[
P_{ij}''(t)+s^2\,P_{ij}(t)=0
\]
for all $j\in \{2,\cdots,n-1\}$. Therefore, we have
\[
P_{i1}(1) = \delta_{i1} + (D^2_{S^{n-1}}\bar u)(e_i,e_1),
\]
and
\[
P_{ij}(1) = \cos(s)\,\delta_{ij} + \frac{\sin (s)}{s} (D^2_{S^{n-1}}\bar u)(e_i,e_j)
\]
for $j\in \{2,\cdots,n-1\}$. On the other hand, we have
\begin{align*}
&D^2_{S^{n-1}}\bar u(e_i,e_1) +s\,\frac{\cos (s)}{\sin(s)}\,g_{ S^{n-1}}(e_i,e_1) +\Big(\frac{1}{s^2}-\frac{\cos(s)}{s\,\sin(s)}\Big)\, d\bar u(e_i)\cdot d\bar u(e_1)\\
&=D^2_{S^{n-1}}\bar u(e_i,e_1)+\delta_{i1},
\end{align*}
and
\begin{align*}
&D^2_{S^{n-1}}\bar u(e_i,e_j) +s\,\frac{\cos (s)}{\sin(s)}\,g_{ S^{n-1}}(e_i,e_j) +\Big(\frac{1}{s^2}-\frac{\cos(s)}{s\,\sin(s)}\Big)\, d\bar u(e_i)\cdot d\bar u(e_j)\\
&=D^2_{S^{n-1}}\bar u(e_i,e_j)+\frac{s\,\cos(s)}{\sin(s)}\,\delta_{ij}.
\end{align*}
It follows that
\begin{align*}
&\det (P_{ij}(1))\\
&=\Big(\frac{\sin (s)}{s}\Big)^{n-2} \det\Big( D^2_{S^{n-1}}\bar u +s\,\frac{\cos (s)}{\sin(s)}\,g_{ S^{n-1}} +\Big(\frac{1}{s^2}-\frac{\cos(s)}{s\,\sin(s)}\Big)\, d\bar u\otimes d\bar u \Big).
\end{align*}
This completes the proof.
\end{proof}
\begin{lemma}\label{lem:am-gm}
Let $x\in A$, and let $p\in\Sigma$ be the unique point satisfying $ \Lb(p)=x. $ Then
\[
0\leq \detchb(p)\det D\Phi(x) \leq \Big(\frac{\sin(|\gd^{S^{n-1}}\bar u( x)|)}{|\gd^{S^{n-1}}\bar u( x)|}\Big)^{n-2}.
\]
\end{lemma}

\begin{proof}
Fix $x\in A$ and let $p=\Lb^{-1}(x)$. We set $w=\chib^{-1}(\gd^\Sigma u)(p)$ and $s=|\gd^{S^{n-1}}\bar u(x)|$. Moreover, choose an orthonormal frame $\{e_1,\ldots,e_{n-1}\}$ on $\Sigma$ near $p$, and define $ \bar e_i=e_i-\la e_i,T_0\rg x. $ Then it follows that $\{\bar e_1,\ldots,\bar e_{n-1}\}$ is an orthonormal basis of $T_xS^{n-1}$, and $ d\Lb_p(e_i) =\sum_{j=1}^{n-1}\chib(e_i,e_j)\bar e_j. $ Since $u=\bar u\circ\Lb$, we have $ s=|w|$. Moreover, differentiating the identity
\[
\gd^{S^{n-1}}\bar u(\Lb) = w-\la w,T_0\rg\Lb
\]
along $e_i$ and using $ \la\bar D_{e_i}w,x\rg =-\chib(e_i,w), $ we obtain
\begin{align*}
&D^2_{S^{n-1}}\bar u(x) (d\Lb_p(e_i),\bar e_i)\\
&= \la\gd^\Sigma_{e_i}w,e_i\rg +\la e_i,T_0\rg\chib(e_i,w) -\la w,T_0\rg\chib(e_i,e_i).
\end{align*}
Summing over $i$ gives
\begin{align*}
\sum_{i,j=1}^{n-1} \chib(e_i,e_j) D^2_{S^{n-1}}\bar u(\bar e_i,\bar e_j) =\operatorname{div}_\Sigma w +\chib(T_0^{\tan},w) -\trchb\la T_0^{\tan},w\rg.
\end{align*}
Combining this identity with Lemma \ref{lem:determinant} and ~\eqref{eq:aux-PDE}, we have
\begin{align*}
&\sum_{i,j=1}^{n-1}\chib(e_i,e_j) \Big( D^2_{S^{n-1}}\bar u(\bar e_i,\bar e_j) +s\frac{\cos(s)}{\sin(s)}\,g_{S^{n-1}}(\bar e_i, \bar e_j)\notag\\
&\qquad\qquad+ \Big(\frac{1}{s^2}-\frac{\cos(s)}{s\,\sin(s)}\Big) d\bar u(\bar e_i)\cdot d\bar u(\bar e_j) \Big)\notag\\
&=(n-1) -\frac34\big( \trchb|T_0^{\tan}|^2 -\chib(T_0^{\tan},T_0^{\tan}) \big)\notag\\
&\quad- \big( \trchb\la T_0^{\tan},w\rg -\chib(T_0^{\tan},w) \big)\notag\\
&\quad- \Big(\frac{1}{s^2}-\frac{\cos(s)}{s\,\sin(s)}\Big) ( \trchb|w|^2-\chib(w,w) )\notag\\
&\leq(n-1) -\frac13\Big( \trchb \Big|w+\frac32T_0^{\tan}\Big|^2- \chib\Big( w+\frac32T_0^{\tan}, w+\frac32T_0^{\tan} \Big) \Big)\notag\\
&\leq n-1.
\end{align*}
Here, in the last inequality, we are using the fact that $\trchb|V|^2-\chib(V,V)\geq0$ for every $V\in T_p\Sigma$, and $\frac{1}{s^2}-\frac{\cos(s)}{s\,\sin(s)}\geq \frac{1}{3}$ for all $0\leq s<\pi$. Hence, Lemma \ref{lem:determinant} and the arithmetic-geometric mean inequality gives
\begin{align*}
0&\leq \detchb(p)\det D\Phi(x)\\
&\leq \Big(\frac{\sin(s)}{s}\Big)^{n-2}\Big[\frac{1}{n-1} \sum_{i,j=1}^{n-1}\chib(e_i,e_j) \Big( D^2_{S^{n-1}}\bar u(\bar e_i,\bar e_j) +s\,\frac{\cos(s)}{\sin(s)}\,g_{S^{n-1}}(\bar e_i, \bar e_j)\\
&\qquad+ \Big(\frac{1}{s^2}-\frac{\cos(s)}{s\,\sin(s)}\Big) d\bar u(\bar e_i)\cdot d\bar u(\bar e_j) \Big)\Big]^{n-1} \\
&\leq \Big(\frac{\sin(s)}{s}\Big)^{n-2}.
\end{align*}
This completes the proof.
\end{proof}
We now complete the proof of Proposition \ref{prop:aux-inequality}. Lemma \ref{lem:surjective} implies $\Phi(A)=S^{n-1}$, the area formula and Lemma \ref{lem:am-gm} gives
\begin{align*}
|S^{n-1}| &\leq \int_A\det D\Phi(x)\,d\vol_{S^{n-1}}(x)\notag\\
&\leq \int_A\Big(\frac{\sin(|\gd^{S^{n-1}}\bar u( x)|)}{|\gd^{S^{n-1}}\bar u(x)|}\Big)^{n-2}\, \det\big( \chib(\Lb^{-1}(x)) \big)^{-1} \,d\vol_{S^{n-1}}(x)\notag\\
&\leq \int_{S^{n-1}} \det\big( \chib(\Lb^{-1}(x)) \big)^{-1} \,d\vol_{S^{n-1}}(x)\notag\\
&= |\Sigma|.
\end{align*}
Here, in the last identity we are using the fact that the Jacobian determinant of the map $ \Lb:\Sigma\longrightarrow S^{n-1} $ is $\detchb$. On the other hand, our scaling~\eqref{eq:normalization} is given by $ \int_{\Sigma} ( \trchb +\frac34\trchb|T_0^{\tan}|^2- \frac34 {\chib} (T_0^{\tan},T_0^{\tan}) )d\mu = (n-1)|\Sigma|$. Hence, we conclude that
\begin{align*}
&\int_{\Sigma} \Big( \trchb +\frac34\trchb|T_0^{\tan}|^2- \frac34 {\chib} (T_0^{\tan},T_0^{\tan}) \Big)d\mu \\
&= (n-1)|\Sigma|\geq (n-1)|S^{n-1}|^{\frac1{n-1}} |\Sigma|^{\frac{n-2}{n-1}}.
\end{align*}
This proves Proposition \ref{prop:aux-inequality} under the case when $\Sigma$ is connected. If $\Sigma$ is not connected, then Proposition follows from applying the inequality to each individual connected component of $\Sigma$ and then summing over all connected component. Since $a^{\frac{n-2}{n-1}}+b^{\frac{n-2}{n-1}}>(a+b)^{\frac{n-2}{n-1}}$ for $a,b>0$, we conclude that
\[
\int_{\Sigma} \Big( \trchb +\frac34\trchb|T_0^{\tan}|^2- \frac34 {\chib} (T_0^{\tan},T_0^{\tan}) \Big)d\mu>(n-1)|S^{n-1}|^{\frac1{n-1}} |\Sigma|^{\frac{n-2}{n-1}}
\]
if $\Sigma$ is disconnected. This completes the proof.

\end{document}